\documentclass[10pt]{article}
\usepackage{amsmath,amssymb,amscd,verbatim,amsthm}
\usepackage{amsmath}
\usepackage{amssymb}

\begin{document}

\begin{center}
{\Large \textbf{On characterizations of $\mathcal{MT}(\lambda)$-functions}}\\%
[0.2in]
\textbf{Wei-Shih Du}\footnote{{\small E-mail address:
wsdu@nknucc.nknu.edu.tw, wsdu@mail.nknu.edu.tw; Tel: +886-7-7172930 ext
6809; Fax: +886-7-6051061.}}

\bigskip

{Department of Mathematics, National Kaohsiung Normal University, Kaohsiung
82444, Taiwan}\bigskip
\end{center}

\hrule\vspace{0.2cm}\bigskip

\noindent \textbf{Abstract:} In this paper, we introduce and share the new
concept of $\mathcal{MT}(\lambda )$-functions and its some
characterizations.\medskip

\medskip

\vspace{0.05cm} \hrule\vspace{0.1cm}\bigskip \bigskip

\noindent {\large \textbf{1. Introduction and preliminaries}}\bigskip

Let $f$ be a real-valued function defined on $\mathbb{R}$. For $c\in \mathbb{%
R}$, we recall
\begin{equation*}
\limsup_{x\rightarrow c^{+}}f(x)=\inf_{\varepsilon
>0}\sup_{c<x<c+\varepsilon }f(x)\text{.}
\end{equation*}

\noindent \textbf{Definition 1.1 (see [1-6]).}\quad A function $\varphi :$ $%
[0,\infty )\rightarrow $ $[0,1)$ is said to be an $\mathcal{MT}$-$function$
(or $\mathcal{R}$-$function$) if $\limsup\limits_{s\rightarrow t^{+}}\varphi
(s)<1$ for all $t\in \lbrack 0,\infty )$. \bigskip

It is obvious that if $\varphi :$ $[0,\infty )\rightarrow $ $[0,1)$ is a
nondecreasing function or a nonincreasing function, then $\varphi $ is an $%
\mathcal{MT}$-function. So the set of $\mathcal{MT}$-functions is a rich
class. However, it is worth to note that there exist functions which are not
$\mathcal{MT}$-functions. \bigskip

\noindent \textbf{Example 1.2 (see [2]).}\quad Let $\varphi :$ $[0,\infty
)\rightarrow $ $[0,1)$ be defined by
\begin{equation*}
\quad \varphi (t):=\left\{
\begin{array}{cc}
\frac{\sin t}{t} & \text{, if }t\in (0,\frac{\pi }{2}]\  \\
0 & \text{, otherwise}.%
\end{array}%
\right.
\end{equation*}%
Since $\limsup\limits_{s\rightarrow 0^{^{+}}}\varphi (s)=1,$ $\varphi $ is
not an $\mathcal{MT}$-function.

Some characterizations of $\mathcal{MT}$-functions were established by Du
[2, Theorem 2.1] and were applied in fixed point theory; for more detail,
one can refer to [1-6] and references therein. In this paper, we introduce
and share Du's recent work on the new concept of $\mathcal{MT}(\lambda )$%
-functions and its some characterizations.\bigskip \bigskip \bigskip

\noindent {\large \textbf{2. New concept of $\mathcal{MT}(\lambda )$%
-functions and its characterizations}}\bigskip

Recall that a real sequence $\{a_{n}\}_{n\in
\mathbb{N}
}$ is called

\begin{enumerate}
\item[(i)] \textit{eventually strictly decreasing }if there exists $\ell \in
\mathbb{N}
$ such that $a_{n+1}<a_{n}$ for all $n\in
\mathbb{N}
$ with $n\geq \ell $;

\item[(ii)] \textit{eventually strictly increasing }if there exists $\ell
\in
\mathbb{N}
$ such that $a_{n+1}>a_{n}$ for all $n\in
\mathbb{N}
$ with $n\geq \ell $;

\item[(iii)] \textit{eventually nonincreasing }if there exists $\ell \in
\mathbb{N}
$ such that $a_{n+1}\leq a_{n}$ for all $n\in
\mathbb{N}
$ with $n\geq \ell $;

\item[(iv)] \textit{eventually nondecreasing }if there exists $\ell \in
\mathbb{N}
$ such that $a_{n+1}\geq a_{n}$ for all $n\in
\mathbb{N}
$ with $n\geq \ell $.
\end{enumerate}

Very recently, Du [7] proved some new characterizations of $\mathcal{MT}$%
-functions as follows. We give the proof for the sake of completeness and
for the readers' convenience.\bigskip

\noindent \textbf{Theorem 2.1 (see [7]).}\quad \textit{Let }$\varphi :$%
\textit{\ }$[0,\infty )\rightarrow $\textit{\ }$[0,1)$\textit{\ be a
function. Then the following statements are equivalent.}

\begin{enumerate}
\item[(a)] $\varphi $\textit{\ is an }$\mathcal{MT}$\textit{-function.}

\item[(b)] \textit{For each }$t\in \lbrack 0,\infty )$\textit{, there exist }%
$r_{t}^{(1)}\in \lbrack 0,1)$\textit{\ and }$\varepsilon _{t}^{(1)}>0$%
\textit{\ such that }$\varphi (s)\leq r_{t}^{(1)}$\textit{\ for all }$s\in
(t,t+\varepsilon _{t}^{(1)})$\textit{.}

\item[(c)] \textit{For each }$t\in \lbrack 0,\infty )$\textit{, there exist }%
$r_{t}^{(2)}\in \lbrack 0,1)$\textit{\ and }$\varepsilon _{t}^{(2)}>0$%
\textit{\ such that }$\varphi (s)\leq r_{t}^{(2)}$\textit{\ for all }$s\in
\lbrack t,t+\varepsilon _{t}^{(2)}]$\textit{.}

\item[(d)] \textit{For each }$t\in \lbrack 0,\infty )$\textit{, there exist }%
$r_{t}^{(3)}\in \lbrack 0,1)$\textit{\ and }$\varepsilon _{t}^{(3)}>0$%
\textit{\ such that }$\varphi (s)\leq r_{t}^{(3)}$\textit{\ for all }$s\in
(t,t+\varepsilon _{t}^{(3)}]$\textit{.}

\item[(e)] \textit{For each }$t\in \lbrack 0,\infty )$\textit{, there exist }%
$r_{t}^{(4)}\in \lbrack 0,1)$\textit{\ and }$\varepsilon _{t}^{(4)}>0$%
\textit{\ such that }$\varphi (s)\leq r_{t}^{(4)}$\textit{\ for all }$s\in
\lbrack t,t+\varepsilon _{t}^{(4)})$\textit{.}

\item[(f)] \textit{For any nonincreasing sequence }$\{x_{n}\}_{n\in
\mathbb{N}
}$\textit{\ in }$[0,\infty )$\textit{, we have }$0\leq \sup\limits_{n\in
\mathbb{N}
}\varphi (x_{n})<1$\textit{.}

\item[(g)] $\varphi $\textit{\ is a function of contractive factor; that is,
for any strictly decreasing sequence }$\{x_{n}\}_{n\in
\mathbb{N}
}$\textit{\ in }$[0,\infty )$\textit{, we have }$0\leq \sup\limits_{n\in
\mathbb{N}
}\varphi (x_{n})<1$\textit{.}

\item[(h)] \textit{For any eventually nonincreasing sequence }$%
\{x_{n}\}_{n\in
\mathbb{N}
}$\textit{\ in }$[0,\infty )$\textit{, we have }$0\leq \sup\limits_{n\in
\mathbb{N}
}\varphi (x_{n})<1$\textit{.}

\item[(i)] \textit{For any eventually strictly decreasing sequence }$%
\{x_{n}\}_{n\in
\mathbb{N}
}$\textit{\ in }$[0,\infty )$\textit{, we have }$0\leq \sup\limits_{n\in
\mathbb{N}
}\varphi (x_{n})<1$\textit{.}\medskip
\end{enumerate}

\noindent \textbf{Proof.}\quad The equivalence of statements (a)-(g) was
indeed proved in [2, Theorem 2.1]. The implications "(h) $\Rightarrow $ (f)"
and "(i) $\Rightarrow $ (g)" are obvious. Let us prove "(f) $\Rightarrow $
(h)". Suppose that (f) holds. Let $\{x_{n}\}_{n\in
\mathbb{N}
}$ be an eventually nonincreasing sequence in $[0,\infty )$. Then there
exists $\ell \in
\mathbb{N}
$ such that $x_{n+1}\leq x_{n}$ for all $n\in
\mathbb{N}
$ with $n\geq \ell $. Put $y_{n}=x_{n+\ell -1}$ for $n\in
\mathbb{N}
$. So $\{y_{n}\}_{n\in
\mathbb{N}
}$ is a nonincreasing sequence in $[0,\infty )$. By (f), we obtain
\begin{equation*}
0\leq \sup\limits_{n\in
\mathbb{N}
}\varphi (y_{n})<1\text{.}
\end{equation*}%
Let $\gamma :=\sup\limits_{n\in
\mathbb{N}
}\varphi (y_{n})$. Then
\begin{equation*}
0\leq \varphi (x_{n+\ell -1})=\varphi (y_{n})\leq \gamma <1\text{ \ \ for
all }n\in
\mathbb{N}
\text{.}
\end{equation*}%
Hence we get
\begin{equation*}
0\leq \varphi (x_{n})\leq \gamma <1\text{ \ \ for all }n\in
\mathbb{N}
\text{ with }n\geq \ell \text{.}
\end{equation*}%
Let
\begin{equation*}
\eta :=\max \{\varphi (x_{1}),\varphi (x_{2}),\cdots ,\varphi (x_{\ell
-1}),\gamma \}<1\text{.}
\end{equation*}%
Then $\varphi (x_{n})\leq \eta $ for all $n\in
\mathbb{N}
$. Hence $0\leq \sup\limits_{n\in
\mathbb{N}
}\varphi (x_{n})\leq \eta <1$ and (h) holds. Similarly, we can varify "(g) $%
\Rightarrow $ (i)". Therefore, from above, we prove that the statements%
\textit{\ }(a)-(i) are all logically equivalent. The proof is completed.%
\hspace{\fill}$\Box $\bigskip

In [7], Du first introduced the concept of $\mathcal{MT}(\lambda )$%
-functions.\bigskip

\noindent \textbf{Definition 2.2 (see [7]).}\quad Let $\lambda \in $ $(0,1]$%
. A function $\mu :$ $[0,\infty )\rightarrow $ $[0,\lambda )$ is said to be
an $\mathcal{MT}(\lambda )$-$function$ if $\limsup\limits_{s\rightarrow
t^{+}}\mu (s)<\lambda $ for all $t\in \lbrack 0,\infty )$. \bigskip

Clearly, an $\mathcal{MT}$-function is an $\mathcal{MT}(1)$-function. It is
quite obvious that $\mu $\ is an $\mathcal{MT}(\lambda )$-function if and
only if $\lambda ^{-1}\mu $\ is an $\mathcal{MT}$-function$.$ \bigskip

\noindent \textbf{Remark 2.3 (see [8]).}\quad Recall that a function $%
\varphi :$ $[0,\infty )\rightarrow \left[ 0,\frac{1}{2}\right) $ is said to
be a $\mathcal{P}$-$function$ [8] if $\limsup\limits_{s\rightarrow t^{+}}\mu
(s)<\frac{1}{2}$ for all $t\in \lbrack 0,\infty )$. So \bigskip a $\mathcal{P%
}$-function is obviously an $\mathcal{MT}\left( \frac{1}{2}\right) $%
-function.\bigskip

The following characterizations of $\mathcal{MT}(\lambda)$-functions is an
immediate consequence of Theorem 2.1. \bigskip

\noindent \textbf{Theorem 2.4 (see [7]).}\quad \textit{Let }$\lambda \in $%
\textit{\ }$(0,1]$\textit{\ and let }$\mu :$\textit{\ }$[0,\infty
)\rightarrow $\textit{\ }$[0,\lambda )$\textit{\ be a function. Then the
following statements are equivalent.}

\begin{enumerate}
\item[(1)] $\mu $\textit{\ is an }$\mathcal{MT}(\lambda )$\textit{-function.}

\item[(2)] $\lambda ^{-1}\mu $\textit{\ is an }$\mathcal{MT}$\textit{%
-function.}

\item[(3)] \textit{For each }$t\in \lbrack 0,\infty )$\textit{, there exist }%
$\xi _{t}^{(1)}\in \lbrack 0,\lambda )$\textit{\ and }$\epsilon _{t}^{(1)}>0$%
\textit{\ such that }$\mu (s)\leq \xi _{t}^{(1)}$\textit{\ for all }$s\in
(t,t+\epsilon _{t}^{(1)})$\textit{.}

\item[(4)] \textit{For each }$t\in \lbrack 0,\infty )$\textit{, there exist }%
$\xi _{t}^{(2)}\in \lbrack 0,\lambda )$\textit{\ and }$\epsilon _{t}^{(2)}>0$%
\textit{\ such that }$\mu (s)\leq \xi _{t}^{(2)}$\textit{\ for all }$s\in
\lbrack t,t+\epsilon _{t}^{(2)}]$\textit{.}

\item[(5)] \textit{For each }$t\in \lbrack 0,\infty )$\textit{, there exist }%
$\xi _{t}^{(3)}\in \lbrack 0,\lambda )$\textit{\ and }$\epsilon _{t}^{(3)}>0$%
\textit{\ such that }$\mu (s)\leq \xi _{t}^{(3)}$\textit{\ for all }$s\in
(t,t+\epsilon _{t}^{(3)}]$\textit{.}

\item[(6)] \textit{For each }$t\in \lbrack 0,\infty )$\textit{, there exist }%
$\xi _{t}^{(4)}\in \lbrack 0,\lambda )$\textit{\ and }$\epsilon _{t}^{(4)}>0$%
\textit{\ such that }$\mu (s)\leq \xi _{t}^{(4)}$\textit{\ for all }$s\in
\lbrack t,t+\epsilon _{t}^{(4)})$\textit{.}

\item[(7)] \textit{For any nonincreasing sequence }$\{x_{n}\}_{n\in
\mathbb{N}
}$\textit{\ in }$[0,\infty )$\textit{, we have }$0\leq \sup\limits_{n\in
\mathbb{N}
}\mu (x_{n})<\lambda $\textit{.}

\item[(8)] \textit{For any strictly decreasing sequence }$\{x_{n}\}_{n\in
\mathbb{N}
}$\textit{\ in }$[0,\infty )$\textit{, we have }$0\leq \sup\limits_{n\in
\mathbb{N}
}\mu (x_{n})<\lambda $\textit{.}

\item[(9)] \textit{For any eventually nonincreasing sequence }$%
\{x_{n}\}_{n\in
\mathbb{N}
}$\textit{\ in }$[0,\infty )$\textit{, we have }$0\leq \sup\limits_{n\in
\mathbb{N}
}\mu (x_{n})<\lambda $\textit{.}

\item[(10)] \textit{For any eventually strictly decreasing sequence }$%
\{x_{n}\}_{n\in
\mathbb{N}
}$\textit{\ in }$[0,\infty )$\textit{, we have }$0\leq \sup\limits_{n\in
\mathbb{N}
}\mu (x_{n})<\lambda $\textit{.}
\end{enumerate}

\bigskip

\noindent \textbf{Remark 2.5.}\quad \lbrack 8, Lemma 3.1] is a special case
of Theorem 2.4 for $\lambda =\frac{1}{2}$.

\bigskip \bigskip \bigskip

\noindent {\large \textbf{Acknowledgments}}

This research was supported by Grant No. MOST 103-2115-M-017-001 of the
Ministry of Science and Technology of the Republic of China. \bigskip
\bigskip \bigskip

\noindent {\large \textbf{References}}

\begin{enumerate}
\item[{[1]}] W.-S. Du, Some new results and generalizations in metric fixed
point theory, Nonlinear Anal. 73 (2010) 1439-1446.

\item[{[2]}] W.-S. Du, On coincidence point and fixed point theorems for
nonlinear multivalued maps, Topology and its Applications 159 (2012) 49-56.

\item[{[3]}] W.-S. Du, H. Lakzian, Nonlinear conditions for the existence of
best proximity points, Journal of Inequalities and Applications, 2012,
2012:206.

\item[{[4]}] W.-S. Du, On Caristi type maps and generalized distances with
applications, Abstract and Applied Analysis, 2013, Volume 2013, Article ID
407219, 8 pages, http://dx.doi.org/10.1155/2013/407219.

\item[{[5]}] W.-S. Du, E. Karapinar, A note on Caristi-type cyclic maps:
related results and applications, Fixed Point Theory and Applications, 2013,
2013:344.

\item[{[6]}] W.-S. Du, F. Khojasteh, Y.-N. Chiu, Some generalizations of
Mizoguchi-Takahashi's fixed point theorem with new local constraints, Fixed
Point Theory and Applications, 2014, 2014:31.

\item[{[7]}] W.-S. Du, New existence results of best proximity points and
fixed points for $\mathcal{MT}(\lambda )$-functions, submitted.

\item[{[8]}] H.K. Pathak, R.P. Agarwal, Y.J. Cho, Coincidence and fixed
points for multi-valued mappings and its application to nonconvex integral
inclusions, Journal of Computational and Applied Mathematics 283 (2015)
201-217.
\end{enumerate}

\end{document}